%% file: preprint.tex
\documentclass[11pt]{article}
\input{Preamble}

\title{Integer L-Shaped Method with Non-Supporting \\ No-Good Optimality Cuts}
\author[1]{Benjamin P. Riley}
\author[1]{Prodromos Daoutidis \thanks{Corresponding author (daout001@umn.edu)}}
\author[1]{Qi Zhang \thanks{Primary corresponding author (qizh@umn.edu)}}
\affil[1]{Department of Chemical Engineering and Materials Science, \break University of Minnesota, Minneapolis, MN 55455, USA}
\date{}

\begin{document}
\maketitle
\begin{abstract}
Two-stage stochastic mixed-integer linear programs with mixed-integer recourse arise in many practical applications but are computationally challenging due to their large size and the presence of integer decisions in both stages. The integer L-shaped method with alternating cuts is a widely used decomposition algorithm for these problems, relying on optimality cuts generated using subproblems to iteratively refine the master problem. A key computational bottleneck in this approach is solving the mixed-integer subproblems to optimality in order to generate separating cuts. This work proposes a modification to the integer L-shaped method with alternating cuts to allow for efficient generation of no-good optimality cuts that are separating for the current master problem solutions without being supporting hyperplanes of the feasible region. These separating cuts are derived from subproblems that are terminated before the optimal solution is found or proven to be optimal, reducing the computational effort required for cut generation. Additionally, an updated optimality cut generation function is proposed to account for the various complexities introduced by this early termination strategy. The effectiveness of the proposed method is demonstrated through two case studies on industrially relevant problems from the literature, which illustrate its advantages in handling large-scale instances with complex mixed-integer subproblems. In these cases, the method achieves substantial reductions in solution time or optimality gap compared to the standard integer L-shaped method with alternating cuts, with performance improvements that increase with mixed-integer subproblem size and complexity.
\end{abstract}


\section{Introduction}
In this work, we consider two-stage stochastic mixed-integer linear programs (MILPs) of the following form:
\begin{equation} 
\begin{split}
    \textrm{{minimize}} \quad & c^\top x + \hat{c}^\top z + Q(x) 
\\
\textrm{{subject to}} \quad & A x + \hat{A} z \geq b
\\
& x \in \{0,1\}^n
\\
& z \in \mathbb{R}^{m_1} \times \mathbb{Z}^{m_2},
\end{split}
\label{OP} 
\tag{OP}
\end{equation}
where the expected value of the optimal second-stage cost, $Q(x) := \sum_{s\in \mathcal{S}}{p_s Q_s(x)}$, is defined via the following optimization problems:
\begin{equation} 
\begin{split}
     Q_s(x) :=\textrm{{min}} \quad & h^\top_s y_s
\\
\textrm{{s.t.}} \quad & F_s x + G_s y_s \geq d_s 
\\
& y_s \in \mathbb{R}^{p_1} \times \mathbb{Z}^{p_2}.
\end{split}
\label{secondone}
\end{equation}
Here, the first-stage decision variables consist of binary state variables (i.e., first-stage variables that affect the second-stage cost), $x$, and mixed-integer variables that do not directly affect the second stage, $z$. Additionally, the second-stage decision variables for scenario $s$, $y_s$, are mixed-integer. 

For ease of exposition, we assume that \eqref{OP} is feasible with relatively complete recourse (i.e., problem \eqref{secondone} is feasible for all feasible $x$ in \eqref{OP}). When these models are combined, the following deterministic equivalent of the two-stage stochastic program is obtained:
\begin{equation} 
\begin{split}
    \textrm{{minimize}} \quad & c^\top x + \hat{c}^\top z + \sum_{s\in \mathcal{S}}{p_s h^\top_s y_s} 
\\
\textrm{{subject to}} \quad & A x + \hat{A} z \geq b
\\
&F_s x + G_s y_s \geq d_s  \quad \forall \, s \in \mathcal{S}
\\
& x \in \{0,1\}^n
\\
&  z \in \mathbb{R}^{m_1} \times \mathbb{Z}^{m_2}
\\
& y_s \in \mathbb{R}^{p_1} \times \mathbb{Z}^{p_2} \quad \forall \, s \in \mathcal{S}.
\end{split}
\label{P} 
\tag{P}
\end{equation}

One of the most popular decomposition methods for two-stage stochastic programs is Benders decomposition (also known as the L-shaped method). However, it cannot be directly applied to problems of the form \eqref{P} due to the integrality constraints on some second-stage decisions. As such, various extensions of Benders decomposition have been developed to solve such problems. For instance, alternative branch-and-cut schemes that branch on second-stage variables, $y_s$, have been devised to allow Benders decomposition to be applied to the solution of \eqref{P} \cite{Weninger2019, Wolsey2020}. In \cite{Fakhri2017}, the concept of local cuts and global cuts are introduced and used in a novel branch-and-cut scheme to apply Benders cuts to the solution of \eqref{P}. In \cite{Sherali2002}, concepts from the reformulation-linearization technique and lift-and-project cuts are adapted to sequentially form convex hull representations of the subproblems, making the generated Benders cuts tight despite the presence of second-stage integrality constraints. In \cite{Carøe1998}, the L-shaped method for continuous recourse is adapted to problems with integer recourse using nonlinear cuts derived from general duality theory. 

Several alternative algorithms for the solution of \eqref{P} that do not extend the ideas of Benders decomposition have also been proposed. For example, \cite{Sen2005} applies a disjunctive decomposition scheme to generate valid inequalities from the solution of multiple relaxations of the second-stage problem. These ideas are extended in \cite{Sen2006}, where branch-and-cut schemes are incorporated to solve second-stage problems. In \cite{Schultz1998}, an efficient enumeration scheme is proposed for a subset of problems in \eqref{P} to generate a set of solutions containing the optimal solution.

Another set of applicable algorithms are decomposition schemes that can be applied to \eqref{P} despite not being specifically designed for two-stage stochastic programs. For instance, column generation schemes for mixed-integer optimization problems with complicating constraints have been applied to two-stage stochastic programs by treating non-anticipativity constraints as complicating constraints \cite{Allman2021b}. This concept has been further developed for the solution of multi-stage stochastic mixed-integer nonlinear programs with discrete state variables \cite{Rathi2025}. Another algorithm for the solution of multi-stage stochastic mixed-integer programs with discrete state variables is the stochastic dual dynamic integer programming (SDDiP) algorithm, which applies Lagrangian cuts and strengthened Benders cuts to approximate nonconvex cost-to-go functions \cite{Zou2019}.

The integer L-shaped method of Laporte and Louveaux \cite{Laporte1993} is a widely used and extensively discussed algorithm that is specifically tailored to exploit the structure of \eqref{P}. This algorithm has been further developed in numerous works. For instance, \cite{Li2018b} extends the method to consider convex nonlinear constraints and incorporate Lagrangian cuts. In several other works, problem-specific features are applied to the algorithm to modify branching schemes \cite{Sanci2021} or efficiently solve subproblems to generate cuts \cite{Biesinger2016}. In \cite{Laporte1998}, a problem-specific dynamic programming algorithm for vehicle routing is applied to the integer L-shaped method to solve subproblems efficiently. Similarly, \cite{Hoogendoorn2023} applies the same problem-specific algorithm for solving vehicle routing subproblems to generate additional problem-specific cuts. 

Despite these successful applications of the integer L-shaped method where problem-specific features are used, there is a need for further development of the algorithm and related algorithms, as many industrially-relevant stochastic programs are large and lack convenient problem-specific features. This work aims to further develop the integer L-shaped method to improve algorithm scalability and increase the number of problems the algorithm is well-suited for. Here, we specifically focus on an extension of the \emph{integer L-shaped method with alternating cuts} of Angulo, Ahmed, and Dey \cite{Angulo2016}. In this algorithm, Benders cuts that are violated by the current master problem solution are sought before computationally costly mixed-integer subproblems are solved.

Following this conceptually simple yet effective strategy, we propose a further modification to the algorithm to allow for efficient generation of no-good optimality cuts that are separating for current master problem solutions without being supporting hyperplanes of the feasible region (i.e., tight cuts). This is achieved by terminating mixed-integer subproblems before the optimal solution is found or proven to be optimal. Additionally, an updated optimality cut generation function is proposed to account for the various complexities introduced by terminating subproblems early. Two computational case studies are presented to demonstrate the potential advantages provided by the proposed modification.

The remainder of this paper is organized as follows. In Section \ref{ilsm}, we describe the integer L-shaped method with alternating cuts from the literature. Section \ref{modification} describes the
proposed modification to the algorithm and the corresponding updated optimality cut generation function. Sections \ref{cs1} and \ref{cs2} focus on the case studies on modular, relocatable manufacturing units in a supply chain and the design of a renewables-based fuel and power production network, respectively. Finally, we provide concluding
remarks in Section \ref{conclusion}.

\section{Integer L-Shaped Method with Alternating Cuts} \label{ilsm}
The following is a multicut implementation of the \emph{integer L-shaped method with alternating cuts} \cite{Angulo2016} for two-stage stochastic MILPs with binary state variables and mixed-integer recourse variables. The algorithm generates so-called optimality cuts in parallel to form convex under-approximations of the optimal value functions of the second-stage scenarios, sometimes referred to as recourse functions, $Q_s(x)$.

The algorithm works by solving a reformulation of \eqref{P} referred to as the master problem (\hypertarget{MPtag}{MP}):
\begin{align} 
\textrm{{minimize}} \quad & c^\top x + \hat{c}^\top z+ \sum_{s\in \mathcal{S}}{p_s \eta_s} 
\\
\textrm{{subject to}} \quad & A x+ \hat{A} z \geq b
\\
&e^\top_{sk} x + \eta_s \geq f_{sk}  \quad \forall \, s \in \mathcal{S}, \, k \in \mathcal{K}_s \label{optcuts}
\\
& x \in \{0,1\}^n
\\
&  z \in \mathbb{R}^{m_1} \times \mathbb{Z}^{m_2}
\\
& \eta \in \mathbb{R}^{S}_+,
\label{MP} 
\end{align}
where constraints (\ref{optcuts}) are optimality cuts that restrict $\eta_s$ to have a value greater than or equal to $Q_s(x)$ for any feasible $x$. 

For the reformulation to be exact, each cut must satisfy $ Q_s(x) \geq f_{sk} - e^\top_{sk} x $ for all feasible $x$ (i.e., the cuts underestimate $Q_s(x)$). Additionally, for each $s$ and feasible $x^*$, there must exist a constraint in (\ref{optcuts}) such that $  Q_s(x^*)=f_{sk}-e^\top_{sk} x^*$ (i.e., at least one underestimator must be tight at $x^*$).  As seen commonly in the literature \cite{Angulo2016,Wolsey2020}, the generation of these optimality cuts is incorporated in a branch-and-cut framework for the solution of (\hyperlink{MPtag}{MP}). During the solution of (\hyperlink{MPtag}{MP}), the formulation is relaxed, and (\ref{optcuts}) starts as an empty set of constraints. Optimality cuts are then added iteratively via an optimality cut generation function.
%Additionally, $\mathcal{K}_s$ is the set of cuts $k$ for scenario $s$.
The cut generation function involves solving $S= |\mathcal{S}|$ subproblems derived from the reformulation of \eqref{P}. If $x$ were to take a fixed value, \eqref{P} would be separable into $S$ separate MILPs. These MILPs define the optimal value functions of the second-stage scenarios, $Q_s(x)$, presented above. They will be referred to as \emph{mixed-integer subproblems} in the remainder of this work:
\begin{equation} 
\begin{split}
     Q_s(x) :=\textrm{{min}} \quad & h^\top_s y_s
\\
\textrm{{s.t.}} \quad & F_s x + G_s y_s \geq d_s 
\\
& y_s \in \mathbb{R}^{p_1} \times \mathbb{Z}^{p_2}.
\end{split}
\label{ISPs}
\tag{$\mathrm{MISP}_s(x)$}
\end{equation}

The dual programs of the LP relaxations of the mixed-integer subproblems are referred to as the \emph{dual separation problems}:
\begin{equation} 
\begin{split}
    Q^{\textrm{LP}}_s(x):=\textrm{{max}} \quad & (d_s-F_s x)^\top u
\\
\textrm{{s.t.}} \quad & G_s^\top u = h_s 
\\
& u \in \mathbb{R}^{m}_+,
\end{split}
\label{DSP} 
\tag{$\mathrm{DSP}_s(x)$}
\end{equation}
where $Q^{\textrm{LP}}_s(x) \leq Q_s(x)$.
The first of the two types of optimality cuts that are generated by the algorithm and included in (\ref{optcuts}) are \emph{Benders optimality cuts},
\begin{equation}
    \eta_s \geq (d_s-F_s x)^\top u_t \quad \forall \, t \in \mathcal{T}_s,
    \label{bendersoptcut}
\end{equation}
where $u_t$ are the extreme points of the feasible region for ($\mathrm{DSP}_s$). Benders optimality cuts are generated by solving (\ref{DSP}) and using the optimal extreme point $u^*$ to generate a constraint of the form (\ref{bendersoptcut}). These constraints generally do not provide a tight underestimate of $Q_s(x)$ at any point and therefore are not supporting hyperplanes of the feasible region. Regardless, the constraints can be used to enforce valid lower bounds on $\eta_s$ and improve the formulation. 

The second type of optimality cuts that are generated by the algorithm and included in (\ref{optcuts}) are referred to as \emph{no-good optimality cuts} in this work:
\begin{equation}
\eta_s \geq Q_s(x^*)+(L_s-Q_s(x^*)) \left( \sum_{j:x^*_j=0}{x_j} +\sum_{j:x^*_j=1}{(1-x_j)}\right).
\label{iop}
\end{equation}

Here, $L_s$ is a lower bound on $Q_s(x)$ for all $x$ values that are feasible in (\ref{P}). One valid value for $L_s$ is the optimal objective value of
\begin{equation} 
\begin{split}
     \underset{x,y_s}{\textrm{{minimize}}} \quad & h^\top_s y_s
\\
\textrm{{subject to}} \quad & F_s x + G_s y_s \geq d_s 
\\
& x \in \{0,1\}^n
\\
& y_s \in \mathbb{R}^{p_1} \times \mathbb{Z}^{p_2}.
\end{split}
\label{L-value-calc}
\tag{$\text{LBP}_s$}
\end{equation}

No-good optimality cuts are generated by solving (\ref{ISPs}) with $x$ as the current solution to (\hyperlink{MPtag}{MP}), $x^*$, and using the optimal objective value to generate a separating constraint of the form (\ref{iop}). A constraint is said to be \textit{separating} if the constraint is satisfied by all feasible solutions to (\ref{P}) but is violated by the current solution to the relaxation of (\hyperlink{MPtag}{MP}). Unlike Benders optimality cuts, constraints (\ref{iop}) provide tight underestimates of $Q_s(x)$ at $x^*$ and are therefore supporting hyperplanes of the feasible region of (\hyperlink{MPtag}{MP}).

A general branch-and-cut implementation of the integer L-shaped method is presented in Algorithm \ref{alg:cap}. When used with the standard optimality cut generation function, described in Algorithm \ref{alg:2}, the method is referred to as the integer L-shaped method with alternating cuts \cite{Angulo2016}. The algorithm proceeds like a standard branch-and-cut procedure by solving multiple relaxations of the master problem with branching constraints while establishing lower and upper bounds to prune or fathom nodes from the tree. The algorithm does this via the following changes to the typical branch-and-cut algorithm. First, integer-feasible nodes that are not fathomed by infeasibility or bound are not immediately accepted as a new incumbent solution and fathomed. Instead, separating Benders and no-good optimality cuts are generated by solving (\ref{DSP}) and (\ref{ISPs}) for all scenarios $s$ and then added to (\hyperlink{MPtag}{MP}). A certificate of feasibility for the current solution is obtained only if no separating cuts are generated, at which point the upper bound and incumbent solution are updated, and the node is fathomed. Second, before branching is done at the root node, separating Benders cuts are sought at fractional solutions to the root node. One valid way of interpreting this step is that the LP relaxation of \eqref{P} is first solved via Benders decomposition.

\renewcommand{\algorithmicrequire}{\textbf{Input:}}
\renewcommand{\algorithmicensure}{\textbf{Output:}}
\begin{algorithm}
\caption{General Integer L-shaped method}\label{alg:cap}
\begin{algorithmic}[1]
\Require $c,\hat{c},p,A,\hat{A},b,\mathcal{S},\epsilon$
\Ensure Optimal solution $\tilde{x}$ to \eqref{P} and optimal value $UB$
\State Solve \eqref{L-value-calc} to compute $L_s$ for all $s \in \mathcal{S}$
\State Initialize the branch-and-bound tree for solving (\hyperlink{MPtag}{MP}) by adding the root relaxation to the list of leaf nodes. Set $UB$ and $LB$ to appropriate initial bounds %Initialize $V= \varnothing$, and $V^{\text{LP}}=\varnothing$
\While{$\frac{|UB-LB|}{|UB|} > \epsilon$}
\State Select a node from the list of leaf nodes
\State Solve the relaxation of (\hyperlink{MPtag}{MP}) at the current node
\vspace{0.5ex}\If {the current node is infeasible} 
\State Fathom the node and return to line 4
\EndIf
\vspace{0.5ex}\State Record the optimal solution $x^*$, $z^*$, and $\eta^*$. Following the usual branch and cut procedure, update the lower bound of the tree, $LB$, if necessary
%\State $LB \gets \max{(LB, \,c^\top x^* + \hat{c}^\top z^*+ \sum_{s\in \mathcal{S}}{p_s \eta^*_s})}$
\vspace{0.5ex}\If{the current node is the root relaxation}
\State Call the cut generation function (Algorithm 2) with $x^*$
%\While{$\eta^* \neq Q^{\textrm{LP}}_s(x^*)$}
\If{a cut was generated}
\State Add generated cuts to (\hyperlink{MPtag}{MP}) and return to line 4 
\EndIf
\EndIf
\vspace{0.5ex} \If{$\,c^\top x^* + \hat{c}^\top z^*+ \sum_{s\in \mathcal{S}}{p_s \eta^*_s} > UB$}
\State Fathom the node and return to line 4
\EndIf
\vspace{0.5ex}\If{$(x^*,z^*)$ contains components that violate integrality constraints} 
\State Select one to branch on following the usual branch-and-cut procedure
\State Remove the current node and append the new nodes to the list of leaf nodes
\Else
\State Call the cut generation function (Algorithm 2) with $x^*$ \Comment{$UB, \tilde{x}$ updated}
%\While{$\eta^* \neq Q^{\textrm{LP}}_s(x^*)$}
\If{a cut was generated}
\State Add generated cuts to (\hyperlink{MPtag}{MP}) 
\Else
\State Fathom the node
\EndIf
\EndIf
\EndWhile
\end{algorithmic}
\end{algorithm}

While the typical version of the algorithm described in this section generates Benders cuts at fractional solutions to the root relaxation, the algorithm could instead be configured differently. For instance, branching could occur immediately at the root relaxation and all cut generation could be performed at integer-feasible node solutions. Alternatively, separating Benders cuts could be sought at fractional solutions before branching at all nodes in the tree, rather than just the root node, as described in \cite{Wolsey2020}. This work includes the implementation that appears most commonly in the literature.

The standard optimality cut generation function, described in Algorithm \ref{alg:2}, uses an \textit{alternating cut strategy}. In this strategy, dual separation problems are solved first to generate Benders cuts. If any separating Benders cuts are found, the cut generation function returns those cuts and terminates; otherwise, the algorithm proceeds to solve mixed-integer subproblems to generate no-good optimality cuts. This is done to reduce computational effort, as the separating Benders cuts may be sufficient to advance the solution of the master problem---either by fathoming the node by bound or by raising the lower bound of the node sufficiently so that another leaf node is selected next (in this instance, the original node may be fathomed by bound when the upper bound is subsequently lowered). The performance improvement from this strategy may be large, as each (\ref{DSP}) is an LP, whereas each (\ref{ISPs}) is a potentially NP-hard MILP. Additionally, the optimality cut generation function keeps lists of $x$ values, $V$ and $V^\text{LP}$, for which the subproblems have already been solved to avoid duplicate evaluations.

\renewcommand{\algorithmicrequire}{\textbf{Input:}}
\renewcommand{\algorithmicensure}{\textbf{Output:}}
\begin{algorithm}
\caption{Standard multicut optimality cut generation function}\label{alg:2}
\begin{algorithmic}[1]
\Require Current (\hyperlink{MPtag}{MP}) solution $(x^*,z^*,\eta^*)$, $V$, $V^{\text{LP}}$
\Ensure A (potentially empty) set of separating cuts for current (\hyperlink{MPtag}{MP}) solution: $\mathcal{C}$
\State Initialize $\mathcal{C}=\varnothing$
\vspace{0.5ex}\If{$x^*\in V$}
\State \Return $\mathcal{C}$ \Comment{We know $\eta^*_s \geq Q_s(x^*) \quad \forall \,s \in \mathcal{S}$}
\EndIf
\vspace{0.5ex}\If{$x^* \notin V^\text{LP}$}
\For{$s \in \mathcal{S}$} \Comment{Can be parallelized}
\State Solve $(\text{DSP}_s(x^*))$ and record optimal dual solution vector $u^*$
\If{$(d_s-F_s x^*)^\top u^* > \eta^*_s$}
\State $\mathcal{C} \gets \mathcal{C} \,\cup\, \{\eta_s \geq(d_s-F_s x)^\top u^*\} $ 
\EndIf
\EndFor
\vspace{0.5ex}\State $V^\text{LP} \gets V^\text{LP} \, \cup \, \{x^*\}$
\vspace{0.5ex}\If{$\mathcal{C}\neq \varnothing$}
\State \Return $\mathcal{C}$
\EndIf
\EndIf
\If{$x^* \in \{0,1\}^n$}
\For{$s \in \mathcal{S}$}\Comment{Can be parallelized}
\State Solve $(\text{MISP}_s(x^*))$ and record optimal value $Q_s(x^*)$
\If{$Q_s(x) > \eta^*_s$}
\State $\mathcal{C} \gets \mathcal{C} \,\cup\, \left\{ \eta_s \geq Q_s(x^*)+(L_s-Q_s(x^*)) \left( \sum\limits_{j:x^*_j=0}{x_j} +\sum\limits_{j:x^*_j=1}{(1-x_j)}\right) \right\} $ 
\EndIf
\EndFor
\vspace{0.5ex}\If{$c^\top x^* + \hat{c}^\top z^*+ \sum\limits_{s\in \mathcal{S}}{p_s Q_s(x^*)} < UB$}
\State $UB \gets c^\top x^* + \hat{c}^\top z^*+ \sum\limits_{s\in \mathcal{S}}{p_s Q_s(x^*)}$
\State $\tilde{x} \gets x^*$
\EndIf
\State $V \gets V \, \cup \, \{x^*\}$
\EndIf
\State \Return $\mathcal{C}$

\end{algorithmic}
\end{algorithm}
\section{Integer L-Shaped Method with Non-Supporting No-Good \break Optimality Cuts} \label{modification}
The Benders cuts in the integer L-shaped method with alternating cuts provide a fast avenue for generating separating cuts without solving an MILP. As discussed above, these cuts may be sufficient for the solution of (\hyperlink{MPtag}{MP}), and the algorithm may never generate the corresponding no-good optimality cut at the same $x^*$. As detailed in \cite{Angulo2016}, this modification yields speedups of one order of magnitude for some problem instances from the literature.

Consider instances where the lower bound on $\eta_s$ at $x^*$ needed for the advancement of the algorithm is greater than $Q_s^\text{LP}(x^*)$ (and therefore cannot be enforced by a Benders cut) but significantly less than $Q_s(x^*)$. In these instances, the integer L-shaped method with alternating cuts would not be able to efficiently generate a separating cut and would instead need to generate a cut of form \eqref{iop}. This is particularly costly in the context of two-stage stochastic programming, where subproblems are parallelized on a per-scenario basis. In this context, all parallel workers must wait for the last subproblem to be solved before proceeding to the next step. Thus, a single particularly difficult \eqref{ISPs} can lead to a significantly diminished performance advantage from parallelization. In this section, we propose a conceptually simple modification to the algorithm that allows for the efficient generation of cuts that enforce lower bounds between $Q_s^\text{LP}(x^*)$ and $Q_s(x^*)$.

The new cuts used by the modified algorithm are, again, no-good optimality cuts:
\begin{equation}
\eta_s \geq \underline{Q_s}(x^*)+(L_s-\underline{Q_s}(x^*)) \left( \sum_{j:x^*_j=0}{x_j} +\sum_{j:x^*_j=1}{(1-x_j)}\right).
\label{mycuts}
\end{equation}

Cuts \eqref{mycuts} differ from \eqref{iop} in that they use a general lower bound on $Q_s(x^*)$, $\underline{Q_s}(x^*)$, in place of $Q_s(x^*)$ itself. Unlike \eqref{iop}, cuts \eqref{mycuts}  generally do not provide a tight underestimate of $Q_s(x)$ at any point and therefore are not supporting hyperplanes of the feasible region. In this sense, they serve the same purpose as Benders cuts in the algorithm: enforce valid bounds on $\eta_s$ that can be generated efficiently.

In the proposed algorithm, the lower bounds, $\underline{Q_s}(x^*)$, are generated by solving \eqref{ISPs} to a termination criteria other than proof of optimality. In this work, the termination criteria considered are terminal optimality gap and time limit. These termination criteria are gradually loosened (i.e., the time limit is extended or the terminal optimality gap is lessened) to generate sequentially tighter cuts as needed---if necessary, the supporting hyperplane of \eqref{iop} is eventually obtained. The proposed algorithm uses the same general integer L-shaped method of Algorithm \ref{alg:cap} but uses the new optimality cut generation function Algorithm \ref{alg:3} instead of Algorithm \ref{alg:2}.

\renewcommand{\algorithmicrequire}{\textbf{Input:}}
\renewcommand{\algorithmicensure}{\textbf{Output:}}
\begin{algorithm}
\caption{Optimality cut generation function with non-supporting no-good optimality cuts }\label{alg:3}
\begin{algorithmic}[1]
\Require Current (\hyperlink{MPtag}{MP}) solution $(x^*,z^*,\eta^*)$, $V$, $V^{\text{LP}}$, $(\alpha_k)_{k=1}^{K}$, $\overline{t}^\text{init.}$
\Ensure A (potentially empty) set of separating cuts for current (\hyperlink{MPtag}{MP}) solution: $\mathcal{C}$
\State Initialize $\mathcal{C}=\varnothing$
\vspace{0.5ex}\If{$x^*\in V \land l(x^*)=K$}
\State \Return $\mathcal{C}$ \Comment{We know $\eta^*_s \geq Q_s(x^*) \quad \forall \,s \in \mathcal{S}$}
\EndIf
\vspace{0.5ex}\If{$x^* \notin V^\text{LP}$}
\For{$s \in \mathcal{S}$} \Comment{Can be parallelized}
\State Solve $(\text{DSP}_s(x^*))$ and record optimal solution vector $u^*$
\If{$(d_s-F_s x^*)^\top u^* > \eta^*_s$}
\State $\mathcal{C} \gets \mathcal{C} \,\cup\, \{\eta_s \geq(d_s-F_s x)^\top u^*\} $ 
\EndIf
\EndFor
\vspace{0.5ex}\State $V^\text{LP} \gets V^\text{LP} \, \cup \, \{x^*\}$
\vspace{0.5ex}\If{$\mathcal{C}\neq \varnothing$}
\State \Return $\mathcal{C}$
\EndIf
\EndIf
\If{$x^* \in \{0,1\}^n$}
\If{$x^*\notin V$}
\State $V \gets V \, \cup \, \{x^*\}$ 
\State $l(x^*)\gets0$
\State $\overline{t}(x^*)\gets \overline{t}^\text{init.}$
\EndIf
\While{$l(x^*) < K$}
\For{$s \in \mathcal{S}$}\Comment{Can be parallelized}
\State Solve $(\text{MISP}_s(x^*))$ with terminal optimality gap $\alpha_{l(x^*)+1}$ and time limit $\overline{t}(x^*)$ 
\State Record terminal lower bound $\underline{Q_s}(x^*)$ and terminal upper bound $\overline{Q_s}(x^*)$ if one was found
\If{$\underline{Q_s}(x) > \eta^*_s$}
\State $\mathcal{C} \gets \mathcal{C} \,\cup\, \left\{ \eta_s \geq \underline{Q_s}(x^*)+(L_s-\underline{Q_s}(x^*)) \left( \sum\limits_{j:x^*_j=0}{x_j} +\sum\limits_{j:x^*_j=1}{(1-x_j)}\right) \right\} $ 
\EndIf
\EndFor
\vspace{0.5ex}\If{$\overline{Q_s}(x^*)$ was found for all $s$ and $c^\top x^* + \hat{c}^\top z^*+ \sum\limits_{s\in \mathcal{S}}{p_s \overline{Q_s}(x^*)} < UB$}
\State $UB \gets c^\top x^* + \hat{c}^\top z^*+ \sum\limits_{s\in \mathcal{S}}{p_s \overline{Q_s}(x^*)}$
\State $\tilde{x} \gets x^*$
\EndIf
\If{Any subproblems terminated due to time limit}
\State $\overline{t}(x^*) \gets 2 \, \overline{t}(x^*)$
\State $l(x^*) \gets \min\{l(x^*)+1,K-1\}$ 
\Else
\State $l(x^*) \gets l(x^*)+1$
\EndIf
\If{$\mathcal{C}\neq \varnothing$}
\State \Return $\mathcal{C}$
\EndIf
\EndWhile
\EndIf
\State \Return $\mathcal{C}$
\end{algorithmic}
\end{algorithm}

One major change incorporated into the modified optimality cut generation function is the accounting of previous subproblem solutions. In particular, functions $l: V \to \mathbb{N}$ and $\overline{t}:V\to \mathbb{R}$ are added to keep a record of the number of times (\ref{ISPs}) has been solved and the current time limit of the subproblem at each $x$, respectively. Further, $l(x)$ is used to determine the terminal optimality gap used when solving (\ref{ISPs}). This is done via $(\alpha_k)_{k=1}^{K}$, a finite, decreasing sequence of terminal optimality gaps with length $K$ that ends in $0$. The sequence of optimality gaps is restricted to have these properties so that the algorithm will generate the supporting hyperplanes of \eqref{iop} (corresponding to an optimality gap of 0) in finite time. As such, the finite convergence of the integer L-shaped method is maintained after the proposed modifications.

Every time a set of subproblems is solved to an optimality gap of $\alpha_{l(x^*)+1}$ for $x^*$, $l(x^*)$ is increased by $1$. This causes subsequent iterations of $(\text{MISP}_s(x^*))$ to be solved to a tighter terminal optimality gap. Additionally, if any subproblems are terminated due to the time limit, $\overline{t}(x^*)$ is doubled. If the time limit terminates a subproblem when the terminal optimality gap is $0$, $l(x^*)$ is not incremented so that the subproblem will be executed with a terminal optimality gap of 0 again. While not necessary for the algorithm's function, the incomplete branch-and-cut trees used to solve mixed-integer subproblems to nonzero optimality gaps can be saved and then resumed in subsequent iterations at the same $x^*$ to eliminate duplicated computations. This, however, could be highly memory- or storage-intensive depending on the implementation.

The no-good optimality cuts generated in Algorithm \ref{alg:2} are supporting hyperplanes of the feasible region, meaning that a certificate of feasibility for the current solution is obtained once it is determined that the solution does not violate the cut at $x^*$. Conversely, the non-supporting optimality cuts generated in Algorithm \ref{alg:3} might not be violated by an infeasible solution to (\hyperlink{MPtag}{MP}). To address this, a while loop is added to the optimality cut generation function to continue generating tighter cuts when no separating cuts are generated in an iteration. This loop continues while $l(x^*) < K$, indicating that $(\text{MISP}_s(x^*))$ has not been solved to an optimality gap of $0$, and terminates when a separating cut is found.

One additional consequence of the changes to the algorithm is that new feasible solutions and potential incumbent solutions are not generated from the lower bounds used in the no-good optimality cuts. Rather, upper bounds on $(\text{MISP}_s(x^*))$,  $\overline{Q_s}(x^*)$, are obtained for each $s \in S$ and then combined (in lines 23--25) to construct feasible solutions and potential incumbent solutions to (\hyperlink{MPtag}{MP}).

%{\color{red}I need to provide a comment about flexibility of the algorithm in terms of when you generate Benders cuts, when you generate cuts at all (current optimal solution (maybe don't discuss this at all duh), any integer feasible solution, fractional solution, fractional solutions only when at the root node). Additionally the algorithm is written with an ambiguous cut generating function intentionally because that is what we change. it works with both.}

\section{Case Studies}
In this section, the proposed integer L-shaped method with non-supporting no-good optimality cuts is evaluated using two industrially relevant two-stage stochastic programs. In the first case study, the optimal set of investments into modular, relocatable manufacturing units is sought for a supply chain that faces uncertain future demands. In the second case study, an optimal process design is sought for a renewables-based fuel and power production network that faces uncertain power generation profiles. 

In both case studies, the problem sizes are modulated and the performance of the proposed algorithm (Mod. ILS) is compared to that of the integer L-shaped method with alternating cuts from the literature (ILS) and directly solving the full-space problem. Subproblems are solved in parallel in both the integer L-shaped method from the literature and the proposed modified algorithm. 

\subsection{Case Study 1: Supply Chain with Relocatable Manufacturing Facilities} \label{cs1}
The problem considered in this case study consists of selecting modular, relocatable manufacturing units for a supply chain facing uncertain future demands. The objective is to invest in a set of modular manufacturing units that yields a minimal combined capital and expected operating cost over the modeled planning horizon. This problem is adapted from \cite{Allman2020}, in which it is a deterministic problem. The deterministic equivalent of the two-stage stochastic program is presented below:

\begin{align} 
\textrm{minimize} \quad &  \sum_{s \in S} p_s \Bigg( \sum_{j \in J} \sum_{t \in T} \Bigg( \sum_{i \in I} c_{ijt} d_{its} x_{ijts} + \sum_{j' \in J, m \in M} h_{jj'mt} w_{jj'mts} \Bigg) \nonumber \\ &\qquad \qquad \qquad \quad\quad+\sum_{i \in I, t \in T} p_{it} q_{its} d_{its} \Bigg) + \sum_{m \in M} \sum_{k \in K_1} g_m y_{mk}
\\
\textrm{{subject to}} \quad & y_{mk} \leq y_{m,k-1} \quad \forall \, m \in \mathcal{M}, \, k = 2,\ldots,|\mathcal{K}_m| \label{e1}
\\
&z_{ms} = \sum_{k\in \mathcal{K}_m}{y_{mk}} \quad \forall \, m \in \mathcal{M}, \, s \in \mathcal{S} \label{e2}
\\
&\sum_{j \in \mathcal{J}}{x_{ijts}} + q_{its} = 1 \quad \forall \, i \in \mathcal{I}, \, t \in \mathcal{T}, \, s \in \mathcal{S} \label{e3}
\\
&v_{0mts}= z_{ms}+\sum_{t'=1}^{t}{\left(\sum_{j' \in \mathcal{J}}{w_{j'jmts}} - \sum_{j' \in \mathcal{J}}{w_{jj'mts}}\right)} \quad \forall \, m \in \mathcal{M}, \, t \in \mathcal{T}, \, s \in \mathcal{S}\label{e4}
\\
&v_{jmts}= \sum_{t'=1}^{t}{\left(\sum_{j' \in \mathcal{J}}{w_{j'jmts}} - \sum_{j' \in \mathcal{J}}{w_{jj'mts}}\right)} \quad \forall \, j \in \mathcal{J}, \, m \in \mathcal{M}, \, t \in \mathcal{T}, \, s \in \mathcal{S}\label{e5}
\\
&\sum_{\in \mathcal{I}}{d_{its}x_{ijts}} \leq \sum_{m \in \mathcal{M}}{v_{jmts}u_m} \quad \forall \, j \in \mathcal{J}, \, t \in \mathcal{T}, \, s \in \mathcal{S}\label{e6}
\\
&\sum_{\in \mathcal{I}}{d_{its}x_{ijts}} \leq u^\text{max} \quad \forall \, j \in \mathcal{J}, \, t \in \mathcal{T}, \, s \in \mathcal{S}\label{e7}
\\
& y_{mk} \in \{0,1\} \quad \forall \, m \in \mathcal{M}, \, k \in \mathcal{K}_m\label{e8}
\\
& 0 \leq q_{its} \leq 1 \quad \forall \, i \in \mathcal{I}, \, t \in \mathcal{T}, \, s \in \mathcal{S}\label{e9}
\\
& 0 \leq x_{ijts} \leq 1 \quad \forall \, i \in \mathcal{I}, \, j \in \mathcal{J}, \, t \in \mathcal{T}, \, s \in \mathcal{S}\label{e10}
\\
& z_{jms} \in \mathbb{Z}_+ \quad \forall \, j \in \mathcal{J}, \, m \in \mathcal{M}, \, s \in \mathcal{S}\label{e11}
\\
& v_{jmts} \in \mathbb{Z}_+ \quad \forall \, j \in \mathcal{J}, \, m \in \mathcal{M}, \, t \in \mathcal{T}, \, s \in \mathcal{S}\label{e12}
\\
& w_{jj'mts} \in \mathbb{Z}_+. \quad \forall \, j \in \mathcal{J}, \, j' \in \mathcal{J}, \, m \in \mathcal{M}, \, t \in \mathcal{T}, \, s \in \mathcal{S}.\label{e13}
\end{align}

In this model, $y_{mk}$ is a binary variable indicating whether or not manufacturing unit $k$ of size $m$ is selected---these are the only first-stage variables in the model. Second-stage variables are as follows: $z_{ms}$ is an auxiliary integer decision that corresponds to the number of units of size $m$ that are selected, $x_{ijts}$ is the fraction of demand at site $i$ that is satisfied by production site $j$, $q_{its}$ is the fraction of demand at site $i$ that is not delivered, $v_{jmts}$ is the number of units of size $m$ located at site $j$ at time $t$, and $w_{jj'mts}$ is the number of units of size $m$ that are relocated from production site $j$ to $j'$ at time $t$. The objective function of the model includes terms for the cost of production, the cost of relocating units, a penalty cost for unmet demand, and the capital cost of the manufacturing units. 

Constraints \eqref{e1} are symmetry-breaking constraints for the binary investment decisions. Constraints \eqref{e2} relate the binary investment decisions and the auxiliary $z_{ms}$ variables. Constraints \eqref{e3} account for the amount of demand at site $i$ that is delivered. Constraints \eqref{e4} and \eqref{e5} account for the number of units at each site over time. Finally, constraints \eqref{e6} and \eqref{e7} impose bounds on the amount of production that can occur at site $i$; the number of manufacturing units present imposes one bound, and the other is a limit on the total production that can occur at any site. Readers are directed to \cite{Allman2020} for further discussion of the model.

Several changes were made to the model to make it a two-stage stochastic program that is compatible with the proposed algorithm. First, the demands, $d_{its}$, were made to be random, scenario-dependent parameters. Second, the integer decisions corresponding to the number of modules that were selected were moved to the second stage as auxiliary variables, and the choice to invest in individual units was instead modeled in the first stage via binary $y_{mk}$. Additionally, the model was made to have relatively complete recourse by the incorporation of $q_{its}$ and its associated penalty term in the objective function.

The size of the supply chain network and the number of scenarios are varied to create nine different problem sizes. The network size is chosen because it is a key parameter that affects the size and complexity of the subproblems. Then, five problem instances are randomly generated using data from Table 2 in \cite{Allman2020} (which includes the distribution for $d_{its}$) for each of the nine problem sizes. Each instance of each problem size is then solved directly as the given full-space model, using the integer L-shaped method with alternating cuts from the literature (ILS), and using the proposed modified integer L-shaped method with non-supporting no-good optimality cuts (Mod. ILS). The integer L-shaped methods are implemented with 32 parallel workers for the solution of the subproblems. All optimization instances are solved using Gurobi version 12.0 \cite{GurobiOptimization2021}, and all algorithms are implemented in JuMP \cite{Dunning2017} via the Julia programming language \cite{Bezanson2017}. Additionally, $(\alpha_k)=(0.1,0.01,0)$ is selected as the sequence of optimality gaps, and 600 seconds is selected as the initial mixed-integer subproblem time limit. Computational statistics of the investigation are summarized in Table \ref{t1}.

\begin{table}[htbp]
\caption{Summary of computational statistics for the case study on the design of supply chains with modular manufacturing units.}
%\fontsize{12}{16}\selectfont
%\setlength\tabcolsep{3pt}
{\centering
\begin{tabular}{cccrrcrrcrr}
\toprule
\multicolumn{1}{l}{} &  & \multicolumn{3}{c}{\textbf{Full-space}} & \multicolumn{3}{c}{\textbf{ILS}} & \multicolumn{3}{c}{\textbf{Mod. ILS}} \\ 
\cmidrule(lr){3-5} \cmidrule(lr){6-8} \cmidrule(lr){9-11} 
\textbf{$\bm{|\mathcal{I}|/|\mathcal{J}|}$} & \textbf{$\bm{|\mathcal{S}|}$} & \multicolumn{1}{l}{\textbf{NS}} & \multicolumn{1}{l}{$\overline{\textbf{gap}}$} & \multicolumn{1}{l}{\textbf{time}} & \multicolumn{1}{l}{\textbf{NS}} & \multicolumn{1}{l}{$\overline{\textbf{gap}}$} & \textbf{time}  & \multicolumn{1}{l}{\textbf{NS}} & \multicolumn{1}{l}{$\overline{\textbf{gap}}$} & \textbf{time}  \\ \midrule
\multirow{3}{*}{25/5} & 32 & 0 & - & \multicolumn{1}{r}{853} & 0 & - & 223  & 0 & - & 250  \\
 & 64 & 1 & 0.68 & \multicolumn{1}{r}{2,089} & 0 & - & \multicolumn{1}{r}{409}  & 0 & - & 468\\
 & 96 & 0 & - & 2,794 & 0 & - & 716  & 0 & - & 820  \\ \midrule
\multirow{3}{*}{50/10} & 32 & 2 & 0.97 & 4,492 & 1 & 0.57 & 3,428  & 0 & - &1,900  \\
 & 64 & 4 & 27.05 & 6,357 & 2 & 1.09 & 4,430 & 0 & -  &3,302  \\
 & 96 & 5 & 60.17 & - & 2 & 3.20 & 4,291 & 0 & - & 5,101  \\ \midrule
\multirow{3}{*}{75/15} & 32 & 4 & 28.95 & 6,018 & 5 & 7.14 & - & 4 & 4.43 & 7,094  \\
 & 64 & 5 & 66.98 & - & 5 & 29.59 & - & 4 & 4.46 & 7,561  \\
 & 96 & 5 & 99.96 & - & 5(1) & $9.23^\text{a}$ & - & 5 & 3.45 & -  \\ \bottomrule
\end{tabular}%
\vspace{0.1em}
%\begin{tablenotes}
%\item [\raisebox{-1ex}{\scalebox{1.5}{*}}] \hspace{-0.4em} \small{Numbers in parentheses indicate the number of problem instances where no feasible solution was found.}
%\end{tablenotes}
\vspace{0.7ex}

\small{Note: $\overline{\textbf{gap}}$ is given as a percentage, and \textbf{time} is in seconds.}
\vspace{0.7ex}

}
\small{$^\text{a}$ Numbers in parentheses indicate the number of problem instances where no feasible solution was found. These instances were excluded from the calculation of average optimality gap.}\label{t1}
\end{table}

In Table \ref{t1}, \textbf{NS} stands for "not solved" and indicates the number of problem instances not solved to optimality within the 3-hour time limit. The average optimality gaps of instances that were not solved to optimality are presented as percentages in the {$\overline{\textbf{gap}}$} columns. Further, the average computation times (in seconds) of the instances that were solved to optimality are presented in the \textbf{time} columns.

As can be seen from the table, all algorithms manage to solve most problem instances with the small, 30-node network to optimality. However, the computation time of the integer L-shaped methods is several times shorter than what is required to directly solve the full-space models for these instances. 

For the medium, 60-node network, Gurobi manages to solve 4 of 15 problem instances in full space to optimality while the integer L-shaped method from the literature solves 10 and the proposed method solves all 15. Further, it can be seen that the terminal optimality gap of these instances for the full-space model increases rapidly with increasing scenario count; this is an intuitive result as it does not exploit the decomposable structure of the model.

The advantage of the proposed method over the integer L-shaped method from the literature is most clear from the results for the large, 90-node network. At this network size, the modified integer L-shaped method is able to solve the problem instances to a terminal optimality gap several times less than the standard integer L-shaped method. One potentially misleading statistic is the average gap of 9.23\% for the standard algorithm and the largest problem instances, which is the average of the four instances where the algorithm found a feasible solution and upper bound. If the instance with no feasible solution were treated as having a terminal gap of 100\%, then this statistic would be near 30\%, like the instances with 64 scenarios.

\begin{figure}[h]
    \centering
    \includegraphics[width=\linewidth]{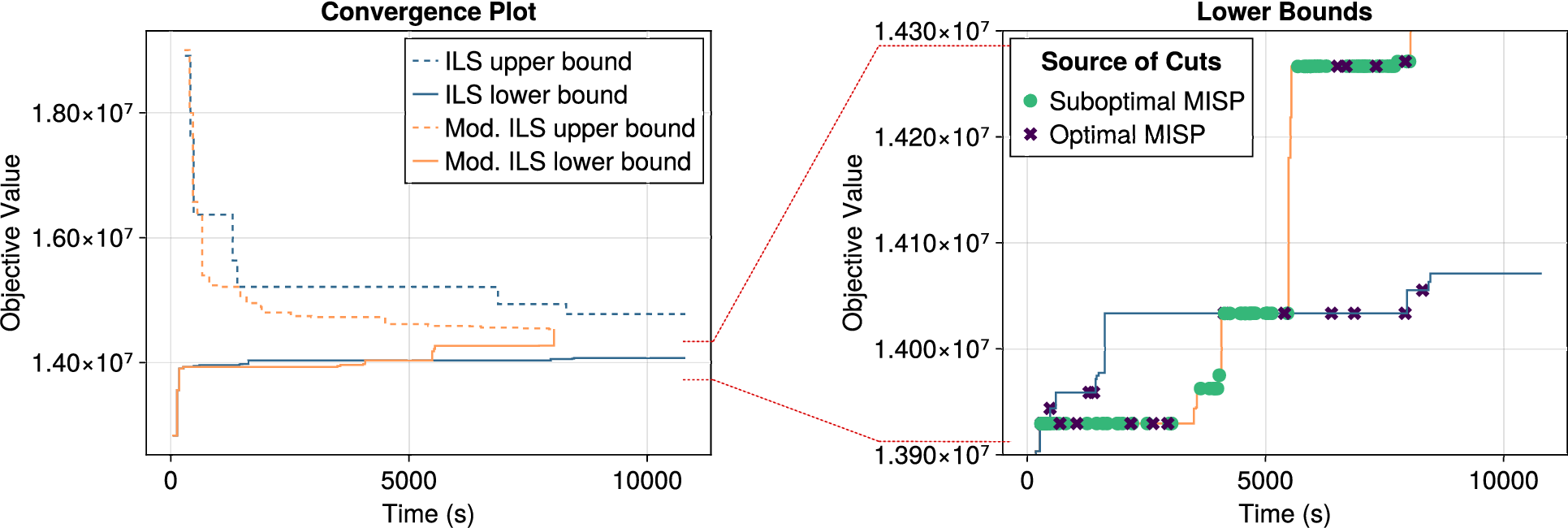}
     \caption{Evolution of upper and lower bounds of the proposed algorithm (Mod. ILS) and standard integer L-shaped method (ILS)}
    \label{f1}
\end{figure}

To evaluate the difference in the behavior of the proposed and standard integer L-shaped methods, the evolution of the upper and lower bounds of the algorithms while solving a specific problem instance with a 90-node network and 64 scenarios is plotted in Figure \ref{f1}.

From the figure, it can be observed that the proposed algorithm solves many more mixed-integer subproblems than the standard algorithm. Additionally, it appears that the lower bound of the standard algorithm remains at the same value for a relatively long time between updates. This demonstrates the key difference between the algorithms: the proposed algorithm generates many relatively loose cuts, while the standard algorithm generates fewer, tighter cuts by solving difficult MILP subproblems to optimality. This takes more than an hour in some subproblem instances. The demonstrated trade-off ultimately provides an advantage that scales with the size and difficulty of the subproblem, as indicated in Table \ref{t1}.

\subsection{Case Study 2: Renewables-Based Fuel and Power Production Network} \label{cs2}
The problem considered in this case study is a superstructure-based network design optimization problem that seeks to design a minimal-cost fuel and power production network that uses intermittent and uncertain wind and solar power. The problem is adapted from \cite{Zhang2019a}, in which it is a deterministic problem. The problem's superstructure models resources (e.g., hydrogen, syngas, and power) and process units as separate nodes. Process units are modeled using multiple operating modes, where some modes correspond to production and others correspond to startup or shutdown. The model includes constraints on the operating modes of each unit, enforcing restrictions such as minimum stay times in some modes, fixed durations of other modes, and predefined transition sequences between some sets of modes. Additionally, storage is available to store resources, and the accumulation of these resources over the planning horizon is modeled. The planning horizon is separated into multiple ``seasons,'' and each season is modeled as a scheduling horizon that cycles multiple times. The deterministic equivalent of the two-stage stochastic program used in this case study is presented in the appendix.

Several changes were made to the model to make it a two-stage stochastic program that is compatible with the proposed algorithm. First, the wind and solar power generation profiles were made to be random, scenario-dependent parameters. Second, storage is assumed to be already built and is therefore not a decision variable. Additionally, first-stage binary decisions are introduced that indicate whether or not process units of fixed nameplate capacities are selected. Furthermore, more operating modes are introduced to the process units (for a total of 5 modes), with each mode having an operating range of 20\% of the unit capacity. Finally, additional constraints on the set of feasible investment decisions are made to improve the formulation of the master problem. These constraints make it so that a process unit cannot be built unless there is a way to produce all of the intermediate resources used by the unit (i.e., it is not possible to build units that cannot be used).
%Readers are referred to \cite{Zhang2019a} for a detailed description of the model

Similar to the first case study, the number of time periods in the scheduling horizon and the number of scenarios are varied to create nine different problem sizes. The number of time periods is chosen because it is a key parameter that affects the size and complexity of the subproblems. Five problem instances are randomly generated from historical wind and solar irradiation data for each of the nine problem sizes. Then, the same three algorithms are used to solve each problem instance. The same choices of $(\alpha_k)$ and $\overline{t}(x)$ are used as in the first case study. A summary of the computational statistics from this case study is presented in Table \ref{t2}.

\begin{table}[htbp]
\caption{Summary of computational statistics for case study on renewables-based fuel and power production network.}
%\fontsize{12}{16}\selectfont
%\setlength\tabcolsep{3pt}
{\centering
\begin{tabular}{cccrrcrrcrr}
\toprule
\multicolumn{1}{l}{} &  & \multicolumn{3}{c}{\textbf{Full-space}} & \multicolumn{3}{c}{\textbf{ILS}} & \multicolumn{3}{c}{\textbf{Mod. ILS}} \\ 
\cmidrule(lr){3-5} \cmidrule(lr){6-8} \cmidrule(lr){9-11} 
\textbf{$\bm{|\mathcal{T}|}$} & \textbf{$\bm{|\mathcal{S}|}$} & \multicolumn{1}{l}{\textbf{NS}} & \multicolumn{1}{l}{$\overline{\textbf{gap}}$} & \multicolumn{1}{l}{\textbf{time}} & \multicolumn{1}{l}{\textbf{NS}} & \multicolumn{1}{l}{$\overline{\textbf{gap}}$} & \textbf{time}  & \multicolumn{1}{l}{\textbf{NS}} & \multicolumn{1}{l}{$\overline{\textbf{gap}}$} & \textbf{time}  \\ \midrule
\multirow{3}{*}{24} & 8 & 2 & 0.61 & \multicolumn{1}{r}{2,963} & 0 & - & 3,496  & 0 & - & 1,641  \\
 & 16 & 5 & 80.60 & - & 0 & - & \multicolumn{1}{r}{9,124}  & 0 & - & 3,105\\
 & 24 & 5 & 99.76 & - & 1 & 8.81 & 5,427  & 0 & - & 3,871  \\ \midrule
\multirow{3}{*}{36} & 8 & 5 & 79.37 & - & 3 & 43.04 & 1,489 & 0 & - &2,094  \\
 & 16 & 5 & 99.99 & -& 3(2) & $10.50^\text{a}$ & 12,148 & 0 & -  &2,418  \\
 & 24 & 5 & 99.87 & - & 4(2) & $97.70^\text{a}$ & 3,915 & 1 & 19.00 & 5,511  \\ \midrule
\multirow{3}{*}{48} & 8 & 5 & 99.99 & - & 2 & 65.94 & 4,807 & 0 & - & 5,850  \\
 & 16 & 5 & 99.99 & - & 3(2) & $61.59^\text{a}$ & 9,535 & 0 & - & 8,289  \\
 & 24 & 5 & 99.99 & - & 2(1) & $96.60^\text{a}$ & - & 0 & - & 10,449 \\ \bottomrule
\end{tabular}
\vspace{0.1em}
%\begin{tablenotes}
%\item [\raisebox{-1ex}{\scalebox{1.5}{*}}] \hspace{-0.4em} \small{Numbers in parentheses indicate the number of problem instances where no feasible solution was found.}
%\end{tablenotes}
\vspace{0.7ex}

\small{\centering Note: $\overline{\textbf{gap}}$ is given as a percentage, and \textbf{time} is in seconds.}
\vspace{0.7ex}

}
\small{$^\text{a}$ Numbers in parentheses indicate the number of problem instances where no feasible solution was found. These instances were excluded from the calculation of average optimality gap.}\label{t2}
\end{table}

Like in the first case study, all algorithms perform well on the smallest problem instances. For all other problem instances, however, Gurobi fails to solve any full-space models to optimality and terminates with an optimality gap of at least 79\%. Clearly, Gurobi does not scale well for the examined problem. The integer L-shaped method from the literature also demonstrates diminished scaling for this problem compared to the first case study, with 7 instances failing to find a feasible solution within the time limit. Additionally, it appears that the performance of the algorithm is highly varied and unpredictable in this case, with some instances being solved to an optimality gap of 10\%, and others of the same size that terminate before a feasible solution is found. This further highlights the potential drawback of the standard integer L-shaped method: the algorithm may spend an excessive amount of time on one $x$ value if there is even one particularly difficult mixed-integer subproblem. This also explains why several instances did not return a feasible solution; the first set of mixed-integer subproblems took excessively long to solve and did not finish before the time limit on the solution of (\hyperlink{MPtag}{MP}) was reached.

Conversely, the proposed modified integer L-shaped method displayed better scaling for this problem than in the first case study. The algorithm was able to solve all but one instance to optimality within the time limit. Additionally, the average computation time increased from 1,641 seconds for problems with 24 time periods and 8 scenarios to 10,449 seconds for problems with 48 time periods and 24 scenarios. This demonstrates the ability of the proposed algorithm to scale well in some problem instances that the standard integer L-shaped method is ill-suited to solve.

\section{Conclusions} \label{conclusion}
In this work, we have presented a modification to the integer L-shaped method with alternating cuts in which mixed-integer subproblems are terminated before the optimal solution is found and proven to be optimal. By doing this, our modified algorithm can efficiently generate no-good optimality cuts that are not supporting hyperplanes of the feasible region, which may be sufficient for fathoming the current node or otherwise advancing the solution of the master problem. This modification is inspired by the alternating cut generation strategy in the literature, which uses Benders cuts to serve the same purpose of efficiently separating the current master problem solution. 

We performed a case study on two industrially relevant two-stage stochastic programs. In the first problem, an optimal set of modular, relocatable manufacturing units is sought for a supply chain network with uncertain future demands. In the second problem, an optimal design is sought for a renewables-based fuel and power production network with uncertain power generation profiles. The computational results of the studies suggest that the modified integer L-shaped method yields better scaling with the size and complexity of the mixed-integer subproblems. Additionally, it appears that the modification yields more consistent performance by avoiding cases where one particularly difficult mixed-integer subproblem halts the progress of the algorithm and mitigates the advantage provided by parallelization.

It should be noted that while this work assumes relatively complete recourse, an extension of the proposed method to instances without relatively complete recourse is conceptually simple via the application of Benders feasibility cuts and no-good cuts. The reader is directed to \cite{Wolsey2020} for details on incorporating these cuts. Additionally, while this work was written from the perspective of two-stage stochastic programming, the algorithm can be applied to any problem that has a suitable structure for the standard integer L-shaped method presented in \cite{Angulo2016}.

\section*{Acknowledgement}
This work was supported by the State of Minnesota through an appropriation from the Renewable Development Account to the University of Minnesota West Central Research and Outreach Center.

\section*{Data availability}
The data and source codes used in Case Study 1 and Case Study 2 are available at \url{https://github.com/ddolab/Modified_Integer_L-Shaped}.
%\section{References}
%Opportunities for future research include a detailed investigation into the selection of optimality gaps, $\alpha$, and time limits $\overline{t}(x)$ used in the algorithm. 

\urlstyle{rm}
\newcommand{\doi}[1]{doi: \url{#1}}
\printbibliography

\newpage\section{Appendix: Model for Case Study 2}
\subsection*{Nomenclature}
\setlength{\LTleft}{0pt}
\emph{Indices / sets}
\begin{longtable}{ll}
    $i\in\mathcal{D}_j$ & set of processes that consume resource $j$ \\
	$h \in H$					& seasons \\
	$i \in I$					& processes \\
	$j \in J$         & resources \\
	$l \in L$         & segments in piecewise-linear approximations \\
	$m \in M$         & operating modes \\
	$t \in T$					& time periods, $T=\{-\theta^{\max}+1,-\theta^{\max}+2,\dots,0,1,\dots,|T|\}$ \\
        $i\in\mathcal{U}_j$ & set of processes that produce resource $j$ 
\end{longtable}
\emph{Subsets}
\begin{longtable}{ll}
	$\bar{J}$						& resources with demand \\
	$\hat{J}$						& resources that must not be discharged \\
	$L_i$               & segments in piecewise-linear approximation for process $i$ \\
	$M_i$								& operating modes in process $i$ \\
	$SQ_i$							& predefined sequences of mode transitions in process $i$ \\
	$\overline{T}_h$		& time periods in season $h$, $\overline{T}_h=\{1,2,\dots,|\overline{T}_h|\}$ \\	
	$TR_i$							& possible mode transitions in process $i$ \\
	$\overline{TR}_{im}$	& modes from which mode $m$ can be directly reached in process $i$ \\
	$\widehat{TR}_{im}$	& modes which can be directly reached from mode $m$ in process $i$ \hspace{5cm}
\end{longtable}
\emph{Parameters}
\begin{longtable}{ll}
  $B^{\max}_{jht}$         & maximum amount of resource $j$ that can be consumed by the process \\
                           & network in time period $t$ of season $h$ \\
  $C^{\max}_i$             & maximum production capacity for process $i$ \\
  $\overline{C}^{\max}_j$  & maximum storage capacity for resource $j$ \\
  $\widehat{C}_{il}$       & production capacity for process $i$ at right end point of segment $l$ \\
  $\widetilde{C}^{\max}_{im}$ & maximum production amount in mode $m$ of process $i$ \\
  $\widetilde{C}^{\min}_{im}$ & minimum production amount in mode $m$ of process $i$ \\
	$D_{jht}$								 & demand for resource $j$ in time period $t$ of season $h$ \\
	$\overline{M}$           & big-M parameter in rate-of-change constraint \\
	$n_h$                    & number of times the representative scheduling horizon of season $h$ is repeated \\ 
	$\alpha_j$              & fixed capital cost for storing resource $j$ \\
	$\delta_{imh}$          & fixed cost for operating in mode $m$ of process $i$ in season $h$ \\
	$\Delta t$							& length of one time period \\	
	$\bar{\Delta}^{\max}_{im}$  & maximum rate of change \\
	$\gamma_{imh}$          & unit cost for operating in mode $m$ of process $i$ in season $h$ \\
	$\epsilon_{jh}$         & fractional loss from storing resource $j$ in season $h$ \\
	$\eta_{ihts}$             & fractional availability of production capacity in process $i$ in time period $t$ of season $h$ \\
	                        & mode $m$ of process $i$ in season $h$ \\
	$\theta_{imm'}$					& minimum stay time in mode $m'$ of process $i$ after switching from mode $m$ to $m'$ \\
	$\bar{\theta}_{imm'm''}$	& fixed stay time in mode $m'$ of the predefined sequence $(m,m',m'')$ in process $i$ \\
	$\theta^{\max}$					& maximum minimum or predefined stay time in a mode \\
	$\rho_{imjh}$           & conversion factor for resource $j$ with respect to the reference resource in \\
	$\phi_{jh}$             & unit cost for purchasing resource $j$ in season $h$ \\
	$\psi_{jh}$             & unit cost for discharging resource $j$ in season $h$
\end{longtable}
\noindent\emph{Unrestricted continuous variables}
\begin{longtable}{ll}
	$\overline{Q}_{jhs}$	& excess inventory for resource $j$ in season $h$
\end{longtable}
\noindent\emph{Nonnegative continuous variables}
\begin{longtable}{ll}
  $B_{jhts}$         & amount of resource $j$ consumed by the process network in time period $t$ of season $h$ \\
	$C_i$							& production capacity for process $i$ \\
    $\tilde{D}_{jhts}$ &  amount of demand of resource $j$ not met in season $h$ and time $t$\\
	$P_{ihts}$         & amount of reference resource produced by process $i$ in time period $t$ of season $h$ \\
	$\overline{P}_{imhts}$  & amount of reference resource produced in mode $m$ of process $i$ in time period $t$ \\
	                       & of season $h$ \\
	$Q_{jhts}$         & inventory level for resource $j$ at time period $t$ of season $h$ \\
	$S_{jhts}$         & amount of resource $j$ discharged from the network in time period $t$ of season $h$ \\
\end{longtable}
\noindent\emph{Binary variables}
\begin{longtable}{ll}
	$x_i$							& equals 1 if process $i$ is selected \\
	$y_{imhts}$				& equals 1 if process $i$ operates in mode $m$ in time period $t$ of season $h$ \\
	$z_{imm'hts}$			& equals 1 if operation of process $i$ switched from mode $m$ to mode $m'$ at time $t$ \\
	                  & of season $h$
\end{longtable}

\newpage\subsection*{Problem Formulation}
\begin{align} 
\textrm{minimize} \quad &  \sum_i \sigma_i  C^\text{nameplate}_i \, x_i+ \sum_s p_s\sum_h \sum_{t \in \overline{T}_h} n_h \Bigg[ \sum_i \sum_{m \in M_i} \left( \delta_{imh} \, y_{imhts} + \gamma_{imh} \, \overline{P}_{imhts} \right)  \nonumber
\\
& \qquad\qquad\qquad\qquad\qquad\qquad\qquad\qquad\,\,\,\,\,\,\,\,\,\,\, \quad\,\,\, + \sum_j \phi_{jh} \, B_{jhts}+\sum_j \psi_{jh} \, S_{jhts} \bigg] 
\\
\textrm{{subject to}} \quad   & C_i = C^\text{nameplate}_i \, x_i \quad \forall \, i \label{a1}\\
  & x_i \in \{0,1\} \quad \forall \, i \label{a2}\\
  & x_i \leq \sum_{i' \in \mathcal{U}_j}{x_{i'}} \quad\forall\, i \in \mathcal{D}_j\label{a3}\\
& \sum_{m \in M_i} y_{imhts} = x_i \quad \forall \, i, \, h, \, t \in \overline{T}_h , \, s\label{a4} \\
  & P_{ihts} = \sum_{m \in M_i} \overline{P}_{imhts} \quad \forall \, i, \, h, \, t \in \overline{T}_h , \, s\label{a5} \\
  & \widetilde{C}^{\min}_{im} \, y_{imhts} \leq \overline{P}_{imhts} \leq \widetilde{C}^{\max}_{im} \, y_{imhts} \quad \forall \, i, \, m \in M_i, \, h, \, t \in \overline{T}_h , \, s\label{a6} \\
  & y_{imhts} \in \{0,1\} \quad \forall \, i, \, m \in M_i, \, h, \, t \in \overline{T}_h, \, s\label{a7}\\
  & -\bar{\Delta}^{\max}_{im} - \overline{M}(2-y_{imhts}-y_{imh,t-1,s}) \leq \overline{P}_{imhts} - \overline{P}_{imh,t-1,s} \nonumber\\
  & \quad \quad \qquad \leq \bar{\Delta}^{\max}_{im} + \overline{M}(2-y_{imhts}-y_{imh,t-1,s}) \quad \forall \, i, \, m \in M_i, \, h, \, t \in \overline{T}_h, \, s\label{a8}\\
  & \sum_{m' \in \overline{TR}_{im}} z_{im'mh,t-1,s} - \sum_{m' \in \widehat{TR}_{im}} z_{imm'h,t-1,s} = y_{imhts} - y_{imh,t-1,s} \quad \forall \, i, \, m, \, h, \, t, \, s \label{a9}\\
  & z_{imm'hts} \in \{0,1\} \quad \forall \, i, \, (m,m') \in TR_i, \, h, \, t \in \overline{T}_h, \, s\label{a10}\\
  & y_{im'hts} \geq \sum_{k=1}^{\theta_{imm'}} z_{imm'h,t-k,s} \quad \forall \, i, \, (m,m') \in TR_i, \, h, \, t \in \overline{T}_h, \, s\label{a11}\\
  &z_{imm'h,t-\bar{\theta}_{imm'm''},s} = z_{im'm''hts} \quad \forall \, i, \, (m,m',m'') \in SQ_i, \, h, \, t \in \overline{T}_h, \, s\label{a12}\\
  & Q_{jhts} = (1-\epsilon_{jh}) Q_{jh,t-1,s} + \sum_i \sum_{m \in M_i} \rho_{imjh} \, \overline{P}_{imhts} + B_{jhts} - S_{jhts} \quad \forall \, j, \, h, \, t, \, s  \label{a13} \\
& P_{ihts} \leq \eta_{ihts} \, C^\text{nameplate}_i \, x_i \quad \forall \, i, \, h, \, t \in \overline{T}_h, \, s \label{a14} \\
& Q_{jhts} \leq \overline{C}_j \quad \forall \, j, \, h, \, t \in \overline{T}_h, \, s \label{a15}\\
& B_{jhts} \leq B^{\max}_{jht} \quad \forall \, j, \, h, \, t \in \overline{T}_h, \, s \label{a16} \\
  & S_{jhts} + \tilde{D}_{jhts}\geq D_{jht} \quad \forall \, j \in \bar{J}, \, h, \, t \in \overline{T}_h, \, s \label{a17}\\
  & S_{jhts} = 0 \quad \forall \, j \in \hat{J}, \, h, \, t \in \overline{T}_h, \, s \label{a18}\\
    & y_{imh,0,s} = y_{imh,|\overline{T}_h|,s} \quad \forall \, i, \, m \in M_i, \, h, \, s \label{a19}\\
  & z_{imm'hts} = z_{imm'h,t+|\overline{T}_h|,s} \quad \forall \, i, \, (m,m') \in TR_i, \, h, \, -\theta^{\max}_i + 1 \leq t \leq -1, \, s\label{a20}\\
  & y_{imh,|\overline{T}_h|,s} = y_{im,h+1,0,s} \quad \forall \, i, \, m \in M_i, \, h \in H \setminus \{|H|\}, \, s \label{a21}\\
  & z_{imm'h,t+|\overline{T}_h|,s} = z_{imm',h+1,ts} \quad \forall \, I, m, m', \, h \in H \setminus \{|H|\}, \, -\theta^{\max}_i + 1 \leq t \leq -1, \, s\label{a22}\\
  & \overline{Q}_{jhs} = Q_{jh,|\overline{T}_h|,s} - Q_{jh,0,s} \quad \forall \, j, \, h, \, s \label{a23} \\
  & Q_{jh,0,s} + n_h \, \overline{Q}_{jhs} = Q_{j,h+1,0,s} \quad \forall \, j, \, h \in H \setminus \{|H|\}, \, s \label{a24}\\
  & Q_{j,|H|,0,s} + n_{|H|} \, \overline{Q}_{j,|H|,s} \geq Q_{j,1,0,s} \quad \forall \, j, \, s \label{a25}
\end{align}

Here, \eqref{a1} are constraints on the investment decisions, $x_i$, that relate the investment decision to the capacity of the process units. Constraints \eqref{a3} are constraints added to the model in this work. They encode the logic that processes should not be built if there is no way to produce an intermediate resource that they consume ($\mathcal{D}_j$ are processes that consume resource $j$, and $\mathcal{U}_j$ are processes that produce $j$). These constraints are added to improve the formulation of the master problem without changing the optimal solution. Constraints \eqref{a4} make it so that a process must be in an operating mode if it is built, and cannot operate if it is not built. Constraints \eqref{a5} model the production of resources to the mode-specific production rate of each process. Constraints \eqref{a6} apply bounds on the mode-specific production rate and ensure that a mode can only produce a resource if the process is operating in that mode. Constraints \eqref{a8} are big-M constraints that apply ramping constraints to the mode-specific production rates. Constraints \eqref{a9} make $z_{imm'hts}$ encode whether or not process $i$ transitions from mode $m$ to mode $m'$ at the end of time period $t$. Constraints \eqref{a11} enforce a minimum amount of time that a process must remain in a mode after transitioning to it. Constraints \eqref{a12} ensure that processes transition between modes at the right time if the modes are part of a fixed sequence of transitions. Constraints \eqref{a13} are mass balances. Constraints \eqref{a14}--\eqref{a18} apply bounds to the various components of the mass balances. Constraints \eqref{a19} and \eqref{a20} enforce the cyclic operation of the system over the courses of each season, $h$. Constraints \eqref{a21} and \eqref{a22} relate the system state at the end of one season to the initial state of the next season. Constraints \eqref{a23}--\eqref{a25} model the accumulation of resources over the course of a cyclic season. The objective function includes components for capital costs, mode-specific operating costs, mode-specific production costs, raw materials costs, and the costs of discharging certain dischargeable resources. Readers are directed to \cite{Zhang2019a} for further discussion of the model components.
\end{document}

%% file: Preamble.tex
  \usepackage[letterpaper, left=1in, right=1in, top=1in, bottom=1in]{geometry}

  \usepackage[affil-it]{authblk}  % for associating author affiliations
  
  \usepackage{palatino}  % font face
  \usepackage{setspace}  % line spacing
  \usepackage{color}  % font color
  \usepackage[dvipsnames]{xcolor}  % adds more color options

	\usepackage[labelfont=bf]{caption}	% caption style

  \usepackage{enumitem}  % control layout of itemize, enumerate, description

  \usepackage{amsmath}  % essential math package
  \usepackage{amssymb}  % math symbols
  \usepackage{amsthm}  % creating proof and theorem environments
  \usepackage{bm}  % boldface letters
  \usepackage{bbm}  % supports \mathbb-style digits, e.g. \mathbbm{1} indicating an all-one vector
  \allowdisplaybreaks[4]  % breaking multiline equations
	\usepackage{relsize}  % can increase the size of math symbols
	\usepackage{algorithm}  % floating environment for algorithmic
	\usepackage{algpseudocodex}  % for easy writing of pseudo code
  
	\theoremstyle{definition}

  \usepackage{pifont}  % dingbats and symbols

	\usepackage{longtable}  % for tables that span more than one page
	\usepackage{multirow}  % for columns spanning multiple rows
	\usepackage{array}  % for column specifications
	\usepackage{booktabs}  % providing nicer tables, using \toprule, \midrule, etc.
    \usepackage{threeparttable}
    
  \usepackage{graphicx}  % enhanced support for graphics, allows e.g. \in­clude­graph­ics
  \usepackage{float}  % better positioning of figures, e.g. by using [h!]
  \usepackage{subfig}  % creating subfigures with \subfloat
  \usepackage{svg}
	\usepackage{hyperref}
	\hypersetup{
		pdfstartview = {FitV},  % fits the width of the page to the window
		colorlinks = true,    % false: boxed links; true: colored links
		linkcolor = NavyBlue,  % color of internal links
		citecolor = NavyBlue,  % color of links to bibliography
		filecolor = NavyBlue,  % color of file links
		urlcolor = NavyBlue  % color of external links
	}

\usepackage[style=numeric-comp, sorting=none, maxbibnames=99, giveninits=true, doi=true, isbn=false, url=false, eprint=true]{biblatex}  % for biblatex, use biber as backend

%% file: Library.bib
@article{Schultz1998,
author = {Schultz, R{\"{u}}diger and Stougie, Leen and van der Vlerk, Maarten H.},
doi = {10.1007/BF02680560},
issn = {0025-5610},
journal = {Math. Program.},
month = jan,
number = {1-3},
pages = {229--252},
title = {{Solving stochastic programs with integer recourse by enumeration: A framework using Gr{\"{o}}bner basis}},
volume = {83},
year = {1998}
}

@article{Sen2005,
author = {Sen, Suvrajeet and Higle, Julia L.},
doi = {10.1007/s10107-004-0566-z},
isbn = {1010700405},
issn = {0025-5610},
journal = {Math. Program.},
month = sep,
number = {1},
pages = {1--20},
title = {{The C3 Theorem and a D2 Algorithm for Large Scale Stochastic Mixed-Integer Programming: Set Convexification}},
volume = {104},
year = {2005}
}

@article{Sen2006,
author = {Sen, Suvrajeet and Sherali, Hanif D.},
doi = {10.1007/s10107-005-0592-5},
issn = {0025-5610},
journal = {Math. Program.},
month = apr,
number = {2},
pages = {203--223},
title = {{Decomposition with branch-and-cut approaches for two-stage stochastic mixed-integer programming}},
volume = {106},
year = {2006}
}

@article{Fakhri2017,
author = {Fakhri, Ashkan and Ghatee, Mehdi and Fragkogios, Antonios and Saharidis, Georgios K.D.},
doi = {10.1016/j.eswa.2017.07.017},
issn = {09574174},
journal = {Expert Syst. Appl.},
month = dec,
pages = {20--30},
publisher = {Elsevier Ltd},
title = {{Benders decomposition with integer subproblem}},
url = {http://dx.doi.org/10.1016/j.eswa.2017.07.017 https://linkinghub.elsevier.com/retrieve/pii/S0957417417304864},
volume = {89},
year = {2017}
}

@inbook{Wolsey2020,
publisher = {John Wiley \& Sons, Ltd},
address ={Hoboken},
author = {Wolsey, Laurence},
isbn = {9781119606475},
title = {Integer Programming},
edition={2},
pages = {235-249},
doi = {https://doi.org/10.1002/9781119606475.ch12},
year = {2020}
}

@article{Sanci2021,
author = {Sanci, Ece and Daskin, Mark S.},
doi = {10.1016/j.trb.2021.01.005},
issn = {01912615},
journal = {Transp. Res. Part B Methodol.},
month = mar,
pages = {152--184},
publisher = {Elsevier Ltd},
title = {{An integer L-shaped algorithm for the integrated location and network restoration problem in disaster relief}},
url = {https://doi.org/10.1016/j.trb.2021.01.005 https://linkinghub.elsevier.com/retrieve/pii/S0191261521000138},
volume = {145},
year = {2021}
}

@article{Hoogendoorn2023,
author = {Hoogendoorn, Y. N. and Spliet, R.},
doi = {10.1287/ijoc.2023.1271},
issn = {15265528},
journal = {INFORMS J. Comput.},
number = {2},
pages = {423--439},
publisher = {http://pubsonline.informs.org/journal/ijoc},
title = {{An Improved Integer L-Shaped Method for the Vehicle Routing Problem with Stochastic Demands}},
volume = {35},
year = {2023}
}

@article{Biesinger2016,
author = {Biesinger, Benjamin and Hu, Bin and Raidl, G{\"{u}}nther},
doi = {10.1016/j.endm.2016.03.033},
issn = {15710653},
journal = {Electron. Notes Discret. Math.},
pages = {245--252},
title = {{An Integer L-shaped Method for the Generalized Vehicle Routing Problem with Stochastic Demands}},
volume = {52},
year = {2016}
}

@incollection{Laporte1998,
address = {Boston, MA},
author = {Laporte, Gilbert and Louveaux, Fran{\c{c}}ois V.},
booktitle = {Fleet Manag. Logist.},
doi = {10.1007/978-1-4615-5755-5_7},
number = {1995},
pages = {159--167},
publisher = {Springer US},
title = {{Solving Stochastic Routing Problems with the Integer L-Shaped Method}},
year = {1998}
}

@article{Carøe1998,
author = {Car{\o}e, Claus C. and Tind, J{\o}rgen},
doi = {10.1007/BF02680570},
issn = {0025-5610},
journal = {Math. Program.},
month = jan,
number = {1-3},
pages = {451--464},
title = {{L-shaped decomposition of two-stage stochastic programs with integer recourse}},
volume = {83},
year = {1998}
}

@article{Sherali2002,
author = {Sherali, Hanif D. and Fraticelli, Barbara M.P.},
doi = {10.1023/a:1013827731218},
issn = {15732916},
journal = {J. Glob. Optim.},
number = {1-4},
pages = {319--342},
title = {{A modification of benders' decomposition algorithm for discrete subproblems: An approach for stochastic programs with integer recourse}},
volume = {22},
year = {2002}
}

@article{Laporte1993,
author = {Laporte, Gilbert and Louveaux, Fran{\c{c}}ois V.},
doi = {10.1016/0167-6377(93)90002-X},
issn = {01676377},
journal = {Oper. Res. Lett.},
month = apr,
number = {3},
pages = {133--142},
title = {{The integer L-shaped method for stochastic integer programs with complete recourse}},
url = {https://linkinghub.elsevier.com/retrieve/pii/016763779390002X},
volume = {13},
year = {1993}
}

@techreport{Weninger2019,
author = {Weninger, Dieter and Wolsey, Laurence},
institution = {Center for Operations Research and Econometrics},
title = {{Benders' algorithm with (mixed)-integer subproblems}},
year = {2019}
}

@article{Li2018b,
author = {Li, Can and Grossmann, Ignacio E.},
doi = {10.1016/j.compchemeng.2018.01.017},
issn = {00981354},
journal = {Comput. Chem. Eng.},
pages = {165--179},
publisher = {Elsevier Ltd},
title = {{An improved L-shaped method for two-stage convex 0–1 mixed integer nonlinear stochastic programs}},
url = {https://doi.org/10.1016/j.compchemeng.2018.01.017},
volume = {112},
year = {2018}
}

@article{Angulo2016,
author = {Angulo, Gustavo and Ahmed, Shabbir and Dey, Santanu S.},
doi = {10.1287/ijoc.2016.0695},
issn = {1091-9856},
journal = {INFORMS J. Comput.},
month = jul,
number = {3},
pages = {483--499},
title = {{Improving the Integer L-Shaped Method}},
url = {https://pubsonline.informs.org/doi/10.1287/ijoc.2016.0695},
volume = {28},
year = {2016}
}

@article{Zou2019,
author = {Zou, Jikai and Ahmed, Shabbir and Sun, Xu Andy},
doi = {10.1007/s10107-018-1249-5},
issn = {14364646},
journal = {Math. Program.},
month = may,
number = {1-2},
pages = {461--502},
publisher = {Springer Verlag},
title = {{Stochastic dual dynamic integer programming}},
url = {https://link-springer-com.ezp1.lib.umn.edu/article/10.1007/s10107-018-1249-5},
volume = {175},
year = {2019}
}

@article{Allman2020,
author = {Allman, Andrew and Zhang, Qi},
doi = {10.1016/J.EJOR.2020.03.045},
issn = {0377-2217},
journal = {Eur. J. Oper. Res.},
month = oct,
number = {2},
pages = {494--507},
publisher = {North-Holland},
title = {{Dynamic location of modular manufacturing facilities with relocation of individual modules}},
volume = {286},
year = {2020}
}

@article{Zhang2019a,
author = {Zhang, Qi and Mart{\'{i}}n, Mariano and Grossmann, Ignacio E.},
doi = {10.1016/J.COMPCHEMENG.2018.06.018},
issn = {0098-1354},
journal = {Comput. Chem. Eng.},
month = mar,
pages = {80--92},
publisher = {Pergamon},
title = {{Integrated design and operation of renewables-based fuels and power production networks}},
volume = {122},
year = {2019}
}

@misc{GurobiOptimization2021,
author = {{Gurobi Optimization LLC}},
title = {{Gurobi Optimizer Reference Manual}},
url = {https://docs.gurobi.com/projects/optimizer/en/current/index.html},
year = {2021}
}

@article{Dunning2017,
archivePrefix = {arXiv},
arxivId = {1508.01982},
author = {Dunning, Iain and Huchette, Joey and Lubin, Miles},
doi = {10.1137/15M1020575},
issn = {0036-1445},
journal = {SIAM Rev.},
month = jan,
number = {2},
pages = {295--320},
publisher = {Society for Industrial and Applied Mathematics Publications},
title = {{JuMP: A Modeling Language for Mathematical Optimization}},
url = {https://epubs.siam.org/doi/10.1137/15M1020575},
volume = {59},
year = {2017}
}

@article{Bezanson2017,
archivePrefix = {arXiv},
arxivId = {1411.1607},
author = {Bezanson, Jeff and Edelman, Alan and Karpinski, Stefan and Shah, Viral B.},
doi = {10.1137/141000671},
issn = {0036-1445},
journal = {SIAM Rev.},
month = jan,
number = {1},
pages = {65--98},
publisher = {Society for Industrial and Applied Mathematics Publications},
title = {{Julia: A Fresh Approach to Numerical Computing}},
url = {https://discourse.julialang.org/t/how-to-cite-julia-in-publications/23060 https://epubs.siam.org/doi/10.1137/141000671},
volume = {59},
year = {2017}
}

@article{Allman2021b,
archivePrefix = {arXiv},
arxivId = {2001.01794},
author = {Allman, Andrew and Zhang, Qi},
doi = {10.1007/s10898-021-01027-w},
issn = {0925-5001},
journal = {J. Glob. Optim.},
month = dec,
number = {4},
pages = {861--880},
publisher = {Springer US},
title = {{Branch-and-price for a class of nonconvex mixed-integer nonlinear programs}},
volume = {81},
year = {2021}
}

@article{Rathi2025,
author = {Rathi, Tushar and Riley, Benjamin P. and Flores-Quiroz, Angela and Zhang, Qi},
doi = {10.1007/s10898-025-01480-x},
issn = {0925-5001},
journal = {Journal of Global Optimization},
month = mar,
pages = {1--32},
publisher = {Springer},
title = {{Column generation for multistage stochastic mixed-integer nonlinear programs with discrete state variables}},
url = {https://link-springer-com.ezp2.lib.umn.edu/article/10.1007/s10898-025-01480-x https://link.springer.com/10.1007/s10898-025-01480-x},
year = {2025}
}
